\documentclass[12pt,leqno]{amsart}
\usepackage{amsmath}
\usepackage{amssymb}
\usepackage{amsthm}
\usepackage{hyperref}
\usepackage{graphicx}
\usepackage{bbm}
\usepackage{setspace}
\usepackage{subcaption}
\usepackage{algorithmic}
\usepackage{algorithm}
\usepackage{bm}

\renewcommand{\arraystretch}{1.2}
\newtheorem{theorem}{Theorem}[section]

\renewcommand{\arraystretch}{1.2}
\def\ba{\begin{array}}
\def\ea{\end{array}}
\def\beq{\begin{equation}}
\def\eeq{\end{equation}}
\def\bea{\begin{eqnarray}}
\def\eea{\end{eqnarray}}
\def\beann{\begin{eqnarray*}}
\def\eeann{\end{eqnarray*}}

\def\R{\mathbb{R}}

\newcommand{\un}[1]{\underline{#1}}
\newcommand{\ov}[1]{\overline{#1}}
\newcommand{\ii}{\textup{i}}

\newcommand{\ignore}[1]{}

\title{\bf Radial Partitioning with Spectral Penalty Parameter Selection in Distributed Optimization for Power Systems}
\author{
Mehdi Karimi }
\date{\today}
 \thanks{\noindent  Mehdi Karimi: Department of Mathematics, Illinois State University, Normal, IL, 61761.  ({e-mail: \bf mkarim3@ilstu.edu}).   Research of the author was supported by	 the US National
Science Foundation (NSF) under Grant No. CMMI-2347120.}

\begin{document}

\begin{abstract}
This paper proposes group-based distributed optimization (DO) algorithms on top of intelligent partitioning for the optimal power flow (OPF) problems. Radial partitioning of the graph of a network is introduced as a systematic way to split a large-scale problem into more tractable sub-problems, which can potentially be solved efficiently with methods such as convex relaxations. Spectral parameter selection is introduced for group-based DO as a crucial hyper-parameter selection step in DO. A software package DiCARP is created, which is implemented in Python using the Pyomo optimization package. Our numerical results for different power network instances show that our designed algorithm returns more accurate solutions to the tested problems with fewer iterations than component-based DO. Our results confirm the importance of smart partitioning and parameter selection for large-scale optimization problems on networks.  
\end{abstract}
\maketitle

\pagestyle{myheadings} \thispagestyle{plain}
\markboth{KARIMI and TUN{\c C}EL}
{Radial Partitioning for Distributed Optimization}

\section{Introduction} \label{sec:intr}
The modern power network is changing with an unprecedented speed by the increasing utilization of new resources of energy such as photovoltaic (PV) generations,  vehicles with vehicle-to-grid capabilities, and battery energy storage.  The problems arising in the operations of power networks are becoming so complex and large that the classical centralized algorithms may  not be applicable anymore. Centralized computation on a single machine for such a complex and expanding problem is intractable. On the other hand, collecting all the information in a centralized controller requires a demanding communication infrastructure and can cause cybersecurity issues.  \emph{Distributed and parallel computation} \cite{bertsekas2015parallel,kshemkalyani2011distributed, lynch1996distributed,wang2021distributed} is a solution in this scenario where the main computational problem is divided into multiple sub-problems, each managed by a computing \emph{agent}. Distributed optimization (DO) algorithms and specifically the Alternating Direction Method of Multipliers (ADMM)\cite{boyd2011distributed} have been extensively studied for power systems, see for example these review papers \cite{molzahn2017survey, maneesha2021survey, wang2017distributed, patari2021distributed} and the references therein. However, despite hundreds of papers in this area, distributed algorithms have not become widely practical. The main issues with distributed optimization algorithms are \emph{slow convergence} and \emph{parameter selection difficulties}, which we consider in this paper.

Not all problems have the potential to be decomposed and treated in a distributed manner. For problems with such a favorable structure, such as many problems in the operations of power networks, the computing agents coordinate and communicate with each other to solve the main problem.
The main optimization problem in the operations of power networks is called Optimal Power Flow (OPF) \cite{frank2012optimal,frank2012optimal2,zhu2015optimization, glover2012power, frank2012primer, skolfield2021operations}, which is a general term for any problem that optimizes an objective function (for example the cost of generation) subject to the physical constraints and engineering limits. OPF covers a wide variety of formulations for the long term and short-term decision-making problems in power system operations.

One of the main contributions of this paper is showing the importance of \emph{partitioning} in the convergence of the DO algorithms. For many problems, there are multiple ways to \emph{partition} them into sub-problems, and the effect of this partitioning on the performance of the distributed algorithms is generally unknown.   Many papers in the literature are \emph{component-based} \cite{mhanna2018component}, where each component of the network is considered a computing agent. In this setup, the number of sub-problems for large scale problems may affect the convergence of the distributed algorithms. Applying a component-based distributed algorithm to a smart grid with thousands of components has major convergence and implementation issues. In contrast, we use the term  \emph{group-based} algorithms that  partition the problem into sub-problems \cite{guo2016case, guo2015intelligent, kyesswa2020impact, murray2021optimal}. The partitioning in  \cite{guo2016case,guo2015intelligent} is based on a heuristic algorithm that uses spectral clustering of graphs and groups together nodes with so called large affinity. \cite{kyesswa2020impact} compares some partitioning techniques from graph theory by numerical experimentation and concludes that they all have fairly similar performances. The authors show that proper partitioning improves the runtime. The same research group shows in a later paper \cite{murray2021optimal} that proper problem partitioning can have a large impact on the convergence rate for the reactive power dispatch problem. The authors in \cite{yang2019distributed} present an ADMM based approach with partitioning for the direct current OPF problem and show by numerical experiments that partitioning can reduce the number of iterations. In this paper, we go beyond partitioning for just reducing the number of sub-problems to propose the idea of \emph{intelligent partitioning}, where we enforce some desirable structure to the sub-problems. One ideal example of intelligent partitioning is splitting a non-convex problem into sub-problems that, even though non-convex, accept exact convex relaxations. This approach can solve sub-problems more efficiently compared to using a generic solver. 
We introduce a systematic \emph{radial partitioning} based on the underlying graph of the network, where each sub-problem induces a tree. This partitioning reduces the number of iterations and creates sub-problems that are not just smaller in size but more tractable compared to the initial problem. 

The second main contribution in this paper is an adaptive parameter selection approach for group-based DO. These parameters heavily affect the convergence of the DO algorithms \cite{mhanna2018component} and different approaches have been proposed such as adaptive \cite{he2000alternating, xu2017adaptive} and machine learning (ML) assisted \cite{zeng2021reinforcement} parameter selection. We propose a simple \emph{Distributed Consensus Algorithm} (DiCA) combined with \emph{Spectral Penalty Parameter Selection} that returns more accurate solutions reported by competent-based algorithms in fewer iterations. The algorithms are implemented in Python using the optimization modeling package Pyomo \cite{hart2011pyomo,bynum2021pyomo}, and are publicly available  as a software package DiCARP \cite{Karimi_DiCARP}. Using this code, we investigate the effect of partitioning on the number of iterations, convergence, and parameter selection for the distributed algorithm.

The main {\bf contributions} of this paper are as follows:
\begin{itemize}
\item Introducing the concept of \emph{intelligent partitioning} for DO and a prototype \emph{radial partitioning}. Designing a graph algorithm for radial partitioning.
\item Designing an adaptive group-based DO for OPF problems, with spectral penalty parameter selection.   
\item Introducing an open-source software package DiCARP for implementing the algorithms in Python, using the Pyomo optimization package. Implementing numerical examples to showcase the performance. 
\end{itemize}

\subsection{Notations} Lower case letters (such as $v$) are used for scalar variables and the bold version of them (such as $\bm v$) are used for a vector of variables. Indexes $i$ and $j$ are used for buses and lines, $k$ for sub-problems, and $t$ for iteration count. For example, $p_{ij}^{k,t}$ is the active power flow on line $ij$ in region $k$ at iteration $t$. As we mostly use the underlying graph of a power network, we use the terms node/bus and edge/line interchangeably. $\ii$ is the imaginary unit, $\Re\{\cdot\}, \Im\{\cdot\}$ are the real and imaginary operators,
$\cdot ^*$ is the complex conjugate operator, and
$\ov \cdot, \un \cdot$ are the upper and lower bound operators.

\section{Group-based Consensus Optimization} \label{sec:group}
In the highest level form, many optimization problems can be written as 
\begin{eqnarray} \label{eq:1}
\min  \sum_{i=1}^N  f_i(x), \ \ \ \ \textup{s.t.}  \ \ \ x \in D,
\end{eqnarray}
where $D \subseteq \R^n$ ($\R^n$ is the Euclidean vector space of dimension $n$) and $f_i: \R^n  \rightarrow \R$ for all $i$. $D$ is the feasible region which is defined by the physical and engineering constraints (for example Kirchhoff's and Ohm's laws). Three important features in the OPF and many other problems in power systems make them very attractive for distributed algorithms: (i) The objective function is highly separated in the sense that in many cases each $f_i$ represents a single component of the network, (ii) For each $i$, only some components of $x$ are involved in the argument of $f_i(x)$, and (iii) The  components of $x$ involved in $f_i(x)$ also appear in a nearly clear-cut subset of constraints defining $D$. Let $A_i$ be the subset of components of $x$ involve in $f_i$ and $\mathcal A_i$ be the linear transformation that takes out the components of $x$ and returns a vector.
Then, by adding some auxiliary variables, the optimization problem can be written as 
\begin{eqnarray} \label{eq:2}
\begin{array}{cl}
\min &  \sum_{i=1}^N  f_i(z_i) \\
 \textup{s.t.}  &  z_i = \mathcal A_i(x), \ \   z_i \in D_i(x), \ \ i\in \{1,\ldots,N\}.
\end{array}
\end{eqnarray}
$D_i(x)$ is the set defined by a subset of constraints that involve $A_i$ and is a function of $x$ since the constraints are not completely decoupled. Equation \eqref{eq:2} is a formulation of our \emph{group-based} DO, where we have $N$ groups and the variables in group $i$ are labeled by the subset $A_i$. There are various approaches to perform this grouping, and the main novelty of our work lies in executing this grouping intelligently. Group-based DO methods contrast with \emph{component-based} DO algorithms, where there is no grouping and each component acts as a computing agent. Given that modern power networks comprise thousands of components, the impact of managing this many sub-problems on the performance of distributed algorithms is not well understood. This framework can be seen as the constrained version of \emph{consensus} optimization discussed in this seminal article by Boyd et al. \cite{boyd2011distributed}, which turns the problem into a form viable for distributed methods. \cite{boyd2011distributed} is a comprehensive article about the ADMM which is the dominant distributed algorithm in the literature. ADMM and many other distributed algorithms are based on the theory of \emph{duality} for optimization and use the \emph{Lagrangian} dual function, which we will introduce later. 


\section{Optimal Power Flow Problem} \label{sec:OPF}
OPF is an optimization problem with an optimization function with a typically linear or quadratic objective function subject to some physical and operational constraints such as Ohm’s Law and Kirchhoff’s Law. Consider the buses and lines of the network as the nodes and edges of a graph $G$, with node set $\mathcal V$ and edge set $\mathcal L$. For each edge set $\mathcal L$, we define $\mathcal L^t$ as the set of all $(j,i)$ such that $(i,j) \in \mathcal L$. Let  $\mathcal G$ be a set of pairs $(g,i)$ where generator $g$ is connected to bus $i$.  The quantities of the network needed to define the OPF problem are given in Table \ref{tab:not}. Assume that $f_i(P_i^g)$ is the objective function associated to $(g,i) \in \mathcal G$, which can be a quadratic generation cost function as $f_i(p_i^g) := c^g_{2,i} (p_i^g)^2 + c^g_{1,i} p_i^g +c^g_{0,i}$. Then we define the optimization problem $OPF(G)$ as the following \cite{mhanna2018component}:
\begin{subequations} \label{eq:opf}
\begin{align}
&\min & &\sum_{(g,i) \in \mathcal G} f_i(p_i^g)  \nonumber \\
& \ \text{s.t.}  & &\un{p}_i^g \leq  p_i^g \leq \ov{p}_i^g; \ \ \   \un{q}_i^g \leq  q_i^g \leq \ov{q}_i^g & (g,i) \in \mathcal G  \label{eq:opf-1}\\
&   & &\un{v}_i \leq  v_i \leq \ov{v}_i, &   i \in \mathcal V  \label{eq:opf-b}\\
&   & &\un{\theta}_{ij} \leq  \theta_i - \theta_j \leq \ov{\theta}_{ij}, &   (i,j) \in \mathcal L \\
&   & & \sum_{(g,i) \in \mathcal G} p_i^g - p_i^d = \sum_{j \in \mathcal V_i} p_{ij} + g_i^{sh} v_i^2,  & i \in \mathcal V \label{eq:opf-c}\\
&   & & \sum_{(g,i) \in \mathcal G} q_i^g - q_i^d = \sum_{j \in \mathcal V_i} q_{ij} - b_i^{sh} v_i^2, & i \in \mathcal V \label{eq:opf-d}\\
&   & & p_{ij} = g_{ij}^c v_i^2-g_{ij} v_i v_j \cos(\Delta \theta _{ij})  + b_{ij} v_i v_j \sin(\Delta \theta _{ij}),   & (i,j) \in \mathcal L  \label{eq:opf-f}\\
&   & & q_{ij} = b_{ij}^c v_i^2-b_{ij} v_i v_j \cos(\Delta \theta _{ij}) - g_{ij} v_i v_j \sin(\Delta \theta _{ij}),   & (i,j) \in \mathcal L \\
&   & & p_{ji} = g_{ji}^c v_i^2-g_{ji} v_j v_i \cos(\Delta \theta _{ji}) + b_{ji} v_j v_i \sin(\Delta \theta _{ji}),   & (i,j) \in \mathcal L \\
&   & & q_{ji} = b_{ji}^c v_i^2-b_{ji} v_j v_i \cos(\Delta \theta _{ji}) + g_{ji} v_j v_i \sin(\Delta \theta _{ji}),   & (i,j) \in \mathcal L \label{eq:opf-j}\\
&   & & \sqrt{p_{ij}^2 + q_{ij}^2} \leq \ov{s}_{ij}; \ \ \ \sqrt{p_{ji}^2 + q_{ji}^2} \leq \ov{s}_{ij},   & (i,j) \in \mathcal L  \label{eq:opf-last}
\end{align}
\end{subequations}
 where $\Delta \theta_{ij} := \theta_i - \theta_j$, $g_{ij}^c = \Re \left\{\frac{Y_{ij}^*-\ii b_{ij}^{ch}/2}{|T_{ij}|^2} \right\}$, $b_{ij}^c = \Im \left\{\frac{Y_{ij}^*-\ii b_{ij}^{ch}/2}{|T_{ij}|^2} \right\}$, $g_{ij} = \Re \left\{\frac{Y_{ij}^*}{T_{ij}} \right\}$, $b_{ij} = \Im \left\{\frac{Y_{ij}^*}{T_{ij}} \right\}$, $g_{ji}^c = \Re \left\{Y_{ji}^*-\ii b_{ji}^{ch}/2 \right\}$, $b_{ji}^c = \Im \left\{Y_{ji}^*-\ii b_{ji}^{ch}/2 \right\}$, $g_{ji} = \Re \left\{\frac{Y_{ji}^*}{T_{ji}^*} \right\}$, and $b_{ji} = \Im \left\{\frac{Y_{ji}^*}{T_{ji}^*} \right\}$. Constraints \eqref{eq:opf-c}-\eqref{eq:opf-d} are power flow equations, \eqref{eq:opf-f}-\eqref{eq:opf-j} are by Kirchhoff’s Law, and \eqref{eq:opf-last} are line flow limits.

 \begin{table}[!t]
\caption{Power flow notations \label{tab:not}}
\centering
\begin{tabular}{|l||l|}
\hline
$Y_{ij}$   & Series admittance (p.u.) of line $ij$ in the $\pi$-model \\ \hline
$v_i \angle \theta_i$ & Complex phasor voltage at bus $i$ \\ \hline
$p_{ij}, q_{ij}$ & Active, reactive power  flow along line $ij$ \\ \hline
$p^g_{i}, q^g_{i}$ &Active, reactive power generation of generator $g$ at bus $i$ \\ \hline
$p^d_{i}, q^d_{i}$ &Active, reactive power demand at bus $i$ \\ \hline
$\bar s_{ij}$ & Apparent power limit of line $ij$ \\ \hline
$\ov{\theta}_{ij}, \un{\theta}_{ij}$ & Upper, lower limits of the difference  \\  & of voltage angles of
buses $i$ and $j$ \\ \hline
$T_{ij}$ & Complex tap ratio of a phase shifting transformer \\ \hline
$b^{sh}_{i}, g^{sh}_{i}$ & Shunt susceptance, conductance at bus $i$ \\ \hline
$b_{ij}^{ch}$ & Charging susceptance of line $ij$ in the $\pi$-model \\ \hline
\end{tabular}
\end{table}

\section{Radial Partitioning} \label{sec:RP}

 The OPF problems are non-convex, difficult to solve, and in general NP-hard \cite{bienstock2019strong,lehmann2015ac}. Several numerical techniques have  been used to solve OPF including Newton method, linear programming, and interior-point methods \cite{zhu2015optimization,mehta2016recent}. Distributed algorithms partition a large-scale OPF problem into smaller sub-problems; however, each of these sub-problems has a similar structure in a smaller size and may still be hard to solve. Most of the papers in the literature need to use a generic solver for a subset of the sub-problems in the distributed algorithm.   

Relaxations and approximations for OPF have been extensively studied (see for example a comprehensive review by \cite{molzahn2019survey}). Convex relaxations such as Second Order Cone Programming (SOCP) and Semidefinite Programming (SDP) \cite{ben2001lectures} are specially attractive since the new advances in convex optimization let us solve many optimization classes efficiently and fast to high accuracy. An important drawback of SDP and in general convex relaxations is that their solutions are not necessarily optimal for the OPF (there is an optimality gap) and may not even be feasible. Several researchers studied the sufficient conditions that SOCP or SDP relaxations are exact for OPF \cite{lavaei2011zero,low2014convex,low2014convex2,bai2008semidefinite,molzahn2014sparsity}.  

For \emph{radial networks}, where the underlying graph is a tree, the exactness of the SCOP and SDP relaxations is proved under some mild conditions \cite{gan2014exact,wei2017optimal,farivar2013branch}. This is the main motivation of the radial decomposition proposed in this paper. By radial partitioning of a general graph with a hard OPF problem, each sub-problem is a radial or close to be a radial network and can be solved more efficiently by using a generic  solver or with convex relaxations. 

Partitioning a graph into trees can be done from different perspectives \cite{chung1978partitions,telle1997algorithms} . We define \emph {radial partitioning} as a partitioning of the node set where each partition induces a tree.  A trivial radial partitioning is putting each node in a separate set. An interesting question for our application is: \emph{How to find a radial partitioning with minimum number of partitions?} As far as we know, this question has not been answered and is an interesting open question. In the following, we propose the greedy Algorithm \ref{alg:part} that uses a brute force approach for the partitioning and works well for our studied problems. 
\begin{algorithm}[h]
\caption{Greedy Radial Partitioning}\label{alg:part}
\begin{algorithmic}
\STATE {\textsc{Input:}} A non-empty connected graph $G$, $k=1$ 
\STATE {\textsc{Start:}} {\bf Define the sub-graph node set $\mathcal V_k=\{  \}$, a empty stack $S=[ \  ]$.} 
\STATE {\bf select a random node $v$ in $G$.}
\STATE $\mathcal V_k \gets \{v\}$.
\STATE {\bf push all neighbors of $v$ into $S$ with $v$ as their parent. }
\STATE {\textsc{While}} ($S$ is non-empty):
\STATE \hspace{0.5cm} {\bf pop an element $u$ from $S$}
\STATE \hspace{0.5cm} IF ({\bf no non-parent neighbor of $u$ is in $\mathcal V_k$}):
\STATE \hspace{1.0cm} $\mathcal V_k \gets \mathcal V_k \cup \{u\}$
\STATE \hspace{0.5cm} {\bf push all neighbors of $u$ not in $\mathcal V_k$ into $S$}
\STATE {\bf save $\mathcal V_k$ as a sub-graph}
\STATE $G \gets G \backslash \mathcal V_k$
\STATE \textsc{If} ({\bf $G$ is non-empty}):
\STATE \hspace{0.5cm} $k \leftarrow k+1$, {\bf go to START}
\STATE \textsc{End}: {\bf return all the saved sub-graphs}
\end{algorithmic}
\end{algorithm}

\begin{figure}[htp]

\subfloat[]{%
  \includegraphics[width=1.5in]{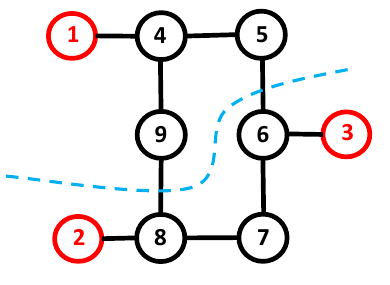}%
}
\qquad
\subfloat[]{%
  \includegraphics[width=3.5in]{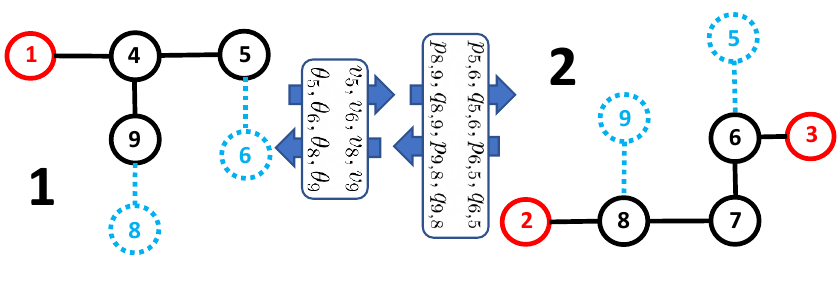}%
}

\caption{(a) Graph of the problem case9 from the MATPOWER library. The red nodes are the buses with generators. The dashed line shows a radial partitioning into two partitions. (b) The node set for the two sub-graphs, both induce trees. The dashed lines and nodes appear in both regions. The variables shown are the ones each region needs to send to the other one for the DiCA.  \label{fig:rad}}
\label{fig:part}
\end{figure}

\begin{theorem}
Algorithm 1 partitions a connected graph into sub-graphs that each induces a tree. 
\end{theorem}
\begin{proof}
For every node added to the set $\mathcal V_k$, its neighbors not already in $\mathcal V_k$ are pushed into $S$. Assume by contradiction that the $\mathcal V_k$ after the WHILE induces a cycle. Clearly, when $\mathcal V_k$ contains one node, the sub-graph is a tree. Consider the first time during the WHILE that $\mathcal V_k$ induces a cycle. This can only happen when the added $u$ has a neighbor in $\mathcal V_k$ other than its parent, which is a contradiction to the IF inside the WHILE.  
\end{proof}
Table \ref{tab:numofreg} shows the number of sub-problems of our greedy algorithm for problems in the MATPOWER library \cite{zimmerman2010matpower}. Figure \ref{fig:rad} shows a radial partitioning for the graph of problem case9.  
\begin{table}[!t]
\caption{Number of regions in radial partitioning for different problem cases \label{tab:numofreg}}
\centering
\begin{tabular}{|c||c|c|c|c|c|c|}
\hline
Case & 9 & 14 & 39 & 89 & 118 & 300 \\
\hline \hline
\# of regions & 2  & 3  & 7 & 10 & 23 & 36\\
\hline
\end{tabular}
\end{table}

\section{Distributed Consensus Algorithm (DiCA)} \label{sec:DiCA}
In this section, we assume that an OPF problem is given, and we have used a partitioning on the graph of the network to get $P$ sub-problems\footnote{Our DiCA algorithm works for any partitioning, not just a radial partitioning.}. Let us define $\mathcal V_k$ as the node set of region $k$. Assume that $N(\mathcal V_k)$ is the union of $\mathcal V_k$ and all the neighbors of the nodes in $\mathcal V_k$.  For a node set $\mathcal V_k$, we also define $G(\mathcal V_k)$ as the subgraph of $G$ induced by the nodes $\mathcal V_k$. For the DiCA algorithm, each agent solves an optimization problem with the same set of constraints as OPF for the graph $\bar G_k:=G(N(\mathcal V_k))$, with a modified objective function. Let us define $\bar {\mathcal V_k}$, $\bar {\mathcal L_k}$, and $\bar {\mathcal G_k}$ as the node, line, and generator sets for $\bar G_k$. For each node or line that appears in more than one $G(N(\mathcal V_k))$, we assign an \emph{auxiliary variable} $\beta$, which acts as the reference for that variable. We also assign a \emph{dual variable} $y^k$ for the variable in each region. For example, in Figure \ref{fig:part}-(b), there are four nodes and four lines that appear in both regions. For node $5$, we define the auxiliary variables $\beta_{v_5}$ and $\beta_{\theta_5}$, and dual variables $y^1_{v_5}$ and $y^1_{\theta_5}$ for region 1 and $y^2_{v_5}$ and $y^2_{\theta_5}$ for region 2. For each node $i$, we define $\mathcal M_i$ as the set of all $k$ such that $\bar G_k$ contains $i$. Similarly, for each line $e=ij$, we define $\mathcal M_e$ as the set of all $k$ such that $\bar G_k$ contains $e$. For each node $i$ such that $|\mathcal M_i| > 1$, each region keeps a copy of the corresponding variables. For example in Figure \ref{fig:part}-(b), for $i=9$, region 1 has a variable $v_9^1$ and region 2 has a variable $v_9^2$. Similarly for every line $e$ with $|\mathcal M_e| > 1$. 

Let $\bm v$, $\bm \theta$, $\bm p$, $\bm q$, and $\bm p^G$ be the vectors that contain $v_i$, $\theta_i$, $p_{ij}$, $q_{ij}$, and $p_i^g$ for a given region, and let $\bm w$ be a vector created by stacking these five vectors together. In other words, $\bm w$ is a vector of all variables in a sub-problem.  Similarly, we define $\bm \beta$ and $\bm y$. Then, the \emph{Lagrangian Dual} for region $k$ is
\begin{align*}
& L_k(\bm w,\bm \beta,\bm y)  := \sum_{(g,i) \in \bar{\mathcal G_k}} f_i(p_i^g)  
+ \sum_{i \in \bar{\mathcal {V}_k}, |\mathcal M_i|>1} \sum_{x\in\{v,\theta\}} y_{x_i} (x_i-\beta_{x_i}) +\frac{\rho^k_{x_i}}{2}(x_i-\beta_{x_i})^2 \nonumber\\
&+ \sum_{ij \in \bar{\mathcal {L}_k} \cup \bar{\mathcal {L}_k}^t, |\mathcal M_{ij}|>1} \sum_{x \in \{p,q\}} y_{x_{ij}} (x_{ij}-\beta_{x_{ij}}) +\frac{\rho^k_{x_{ij}}}{2}(x_{ij}-\beta_{x_{ij}})^2 \nonumber\\
\end{align*}
Algorithm \ref{alg:DiCA} describes the DiCA algorithm. At the first stage, which is the computationally demanding one, we solve OPF with $L_k$ as the objective function for each region to update $\bm w$. At the second and third stages, the algorithm updates $\bm \beta$ and $\bm y$ respectively with simple algebraic operations. The forth stage is for updating the penalty parameters. As can be seen, variables $\beta$ do not have the index $k$ for the region, since these are independent of regions and are the references that the local copies \emph{consent to} be equal to them. 
\begin{algorithm}[h!]
\caption{Distributed Consensus Algorithm (DiCA)}\label{alg:DiCA}
\begin{algorithmic}
\STATE {\textsc{Input:}} OPF problem with graph $G$ partitioned into $P$ sub-graphs. 
\STATE {\textsc{Start:}} { Set $t=0$. Initialize $\bm \beta^0$, $\bm y^{k,0}$s, and $\bm \rho^{k,0}$.} 
\STATE {\textsc{While}} (Stopping criteria are not met):
\STATE \hspace{0.5cm} {{\textsc{For} $k=1,\ldots,P$ solve:}}\[  \begin{array}{ccc}
\bm w^{k,t+1} \leftarrow &\text{argmin}_{\bm w}  & L_k(\bm w,\bm \beta^{t}, \bm y^{k,t}) \\
&\text{s.t.} & \text{\eqref{eq:opf-1}-\eqref{eq:opf-last} for  $\bar{\mathcal V_k}$, $\bar{\mathcal L_k}$, and $\bar{\mathcal G_k}$}
\end{array} \]
\STATE \hspace{0.5cm} {{\textsc{For} nodes with $|\mathcal M_{i}| > 1$ and edges with $|\mathcal M_{ij}| > 1$}:}
\[ \ \ \ \ \  \renewcommand{\arraystretch}{1.3} \begin{array}{l}
\beta^{t+1}_{x_i} \leftarrow \frac{1}{|\mathcal M_i|}\sum_{k \in \mathcal M_i} (x^{k,t+1}_i + \frac{1}{\rho_x} y^{k,t}_{x_i} ), \ \ x \in \{v,\theta\}  \\
\beta^{t+1}_{x_e} \leftarrow \frac{1}{|\mathcal M_e|}\sum_{k \in \mathcal M_e} (x^{k,t+1}_e + \frac{1}{\rho_x} y^{k,t}_{x_e} ), \ \ e=ij,ji, x \in \{p,q\} \\
\end{array} \]
\STATE \hspace{0.5cm} {{\textsc{For} nodes with $|\mathcal M_{i}| > 1$ and edges with $|\mathcal M_{ij}| > 1$}:}
\[ \ \ \ \ \  \renewcommand{\arraystretch}{1.3} \begin{array}{l}
y^{k,t+1}_{x_i} = y^{k,t}_{x_i} + \rho^{k,t}_{x_i} (x^{k,t+1}_i - \beta^{t+1}_{x_i}), \ \ x \in \{v,\theta\} \\
y^{k,t+1}_{x_e} = y^{k,t}_{x_e} + \rho^{k,t}_{x_e} (x^{k,t+1}_e - \beta^{t+1}_{x_e}), \ \  e=ij,ji, x \in \{p,q\}\\
\end{array} \]
\STATE \hspace{0.5cm} {{\textsc{For} nodes with $|\mathcal M_{i}| > 1$ and edges with $|\mathcal M_{ij}| > 1$}, update $\rho^{k,t+1}_{x_i}$ and $\rho^{k,t+1}_{x_e}$.}
\STATE \hspace{0.5cm}  $t \leftarrow t+1$
\end{algorithmic}
\end{algorithm}

This algorithm is a simple generalization of the consensus ADMM described in \cite{boyd2011distributed}-Chapter 7, with promising numerical performance. The structure of the algorithm is simpler than many other distributed optimization algorithms for OPF in the literature, which simplifies crucial parameter tuning. 
For the stopping criteria, each agent checks if its own primal and dual residuals are smaller than a given $\epsilon$. This is more practical than the cases where a centralized unit has to calculate the stopping criteria. For each agent $k$, these residual vectors at iteration $t$ are calculated as follows, which are the generalization of the residuals for the general ADMM approach in \cite{boyd2011distributed}:
\begin{eqnarray*} \label{eq:res}
&&\bm r^{k,t} := \bm w^{k,t} - \bm \beta^{k,t}, \ \ \ \ 
\bm d^{k,t} := \bm{\beta_\rho}^{k,t} - \bm{\beta_\rho}^{k,t-1},
 \\  \nonumber 
&& \bm w^{k,t} := \left [\begin{array}{c}
\mathbf v^{k,t}  \\
\bm \theta^{k,t}  \\
\bm p^{k,t} \\
\bm q^{k,t} 
\end{array} \right]_{|\mathcal M|>1} 
\bm \beta^{k,t} := \left [\begin{array}{c}
 \bm{\beta_ v} ^{t} \\
\bm{\beta_ \theta} ^{t} \\
 \bm{\beta_ p} ^{t} \\
 \bm{\beta_ q} ^{t}
\end{array} \right]_k
\bm{\beta_\rho}^{k,t} := \left [\begin{array}{c}
\rho _v \bm{\beta_ v} ^{t} \\
\rho_\theta \bm{\beta_ \theta} ^{t}\\
\rho_p \bm{\beta_ p} ^{t} \\
\rho_q \bm{\beta_ q} ^{t} 
\end{array} \right]_k
\bm y^{k,t} := \left [\begin{array}{c}
 \bm{y_ v} ^{k,t} \\
\bm{y_ \theta} ^{k,t} \\
 \bm{y_ p} ^{k,t} \\
 \bm{y_ q} ^{k,t}
\end{array} \right].
\end{eqnarray*}
We used the notation $[\cdot]_{_{|\mathcal M|>1}}$ to show that only nodes $i$ and lines $e$ with $|\mathcal M_i| > 1$ and $|\mathcal M_e|>1$ are used in calculating the residual vectors. Also  $\bm \beta^{k,t}$ and $\bm{\beta_\rho}^{k,t}$ are defined to contain variable $\beta$s that are involved in region $k$. Then, the stopping criteria for agent $k$ with tolerance $\epsilon$ are:
\begin{eqnarray} \label{eq:res2} 
\|\bm r^{k,t}\|_2  \leq \epsilon \max\{\|\bm w^{k,t}\|_2, \|\bm \beta^{k,t}\|_2\}, \ \ \
\|\bm d^{k,t}\|_2  \leq \epsilon \|\bm y^{k,t}\|_2
\end{eqnarray}
The DiCA stops when \eqref{eq:res2} is satisfied for all the regions. 

\section{Spectral Penalty Parameter} \label{sec:adap}
Spectral adaptive parameter selection for vanilla ADMM was introduced in \cite{xu2017adaptive} and then extended to consensus ADMM in \cite{xu2017adaptive-cons}. The idea developed in \cite{xu2017adaptive} is based on a classical parameter selection method by Barzilai and Borwein \cite{barzilai1988two} that is applied to the dual problem when we use \emph{Douglas-Rachdord splitting (DRS)} for solving it. We get a parameter selection method for the primal problem since in vanilla ADMM, the DRS for the dual problem is equivalent to applying ADMM to the primal problem. The authors in \cite{mhanna2018adaptive} tailored these results to design an adaptive component-based ADMM  algorithm for OPF. In this section, we present the spectral adaptive parameter selection for group-based algorithms. For every node $i$ where $|\mathcal M_i| > 1$, $\rho_{v_i}$ and $\rho_{\theta_i}$ are updated adaptively. Similarly, for every edge $e$ where $|\mathcal M_e| >1$, $\rho_{p_e}$ and $\rho_{q_e}$ are updated adaptively. Let us define $\bar y^{k,t+1}_{x_i} = \bar y^{k,t}_{x_i} + \rho^{k,t}_{x_i} (x^{k,t+1}_i - \beta^{t+1}_{x_i})$	for a variable $x_i$. Then, we define
\begin{eqnarray} \label{eq:adapt-1}
\alpha^{k,t}_{SD,x_i} = \frac{\sum_{k \in \mathcal M_i} (\Delta \bar y_{x_i}^{k,t})^2}{\sum_{k \in \mathcal M_i} (\Delta \bar y_{x_i}^{k,t})(\Delta x_i^{k,t})}, \ \alpha^{k,t}_{MG,x_i} = \frac{\sum_{k \in \mathcal M_i} (\Delta \bar y_{x_i}^{k,t})(\Delta x_i^{k,t})}{\sum_{k \in \mathcal M_i} (\Delta x_i^{k,t})^2},
\end{eqnarray}
where SD stands for
\emph{steepest descent} and MG for \emph{minimum gradient} \cite{xu2017adaptive}, representing two different methods for calculating  the coefficient of the linear term in approximating the Hessian matrix. 
Similar to \cite{xu2017adaptive-cons}, we use the hybrid of these two estimators as
\begin{eqnarray} \label{eq:adapt-2}
\alpha^{k,t}_{x_i} := \left \{\begin{array}{ll}
\alpha^{k,t}_{MG,x_i}    &  \text{if} \ 2\alpha^{k,t}_{MG,x_i} > \alpha^{k,t}_{SD ,x_i}	 \\
\alpha^{k,t}_{SD,x_i}-\alpha^{k,t}_{MG,x_i}/2  &  \text{o.w.}
\end{array} \right.
\end{eqnarray}
For safeguarding our adaptive updates, we can check that the correlation
\begin{eqnarray} \label{eq:adapt-3}
\alpha^{k,t}_{c,x_i} := \frac{\sum_{k \in \mathcal M_i} (\Delta \bar y_{x_i}^{k,t})(\Delta x_i^{k,t})}{\sqrt{\sum_{k \in \mathcal M_i} (\Delta \bar y_{x_i}^{k,t})^2 \sum_{k \in \mathcal M_i} (\Delta x_i^{k,t})^2}}, 
\end{eqnarray} is bounded away from zero. We also define $\beta^{k,t}_{SD,x_i}$, $\beta^{k,t}_{MG,x_i}$, $\beta^{k,t}_{x_i}$, and $\beta^{k,t}_{c,x_i}$ using the same formulas by replacing $\Delta \bar y_{x_i}^{k,t}$ with $\Delta y_{x_i}^{k,t}$ and $\Delta x_i^{k,t}$ with $\Delta \beta _{x_i}^{t}$.
Using the above estimates and safeguarding, the final adaptive formula is 
\begin{eqnarray} \label{eq:adapt-4}
\hat \rho_{x_i}^{k,t+1} &:=& \left\{\begin{array}{ll}
\sqrt{\alpha^{k,t}_{x_i} \beta^{k,t}_{x_i}}   & \text{if $\alpha^{k,t}_{c,x_i} > \epsilon_c$ and $\beta^{k,t}_{c,x_i} > \epsilon_c$} \\
\alpha^{k,t}_{x_i}    & \text{if $\alpha^{k,t}_{c,x_i} > \epsilon_c$ and $\beta^{k,t}_{c,x_i} \leq \epsilon_c$} \\
\beta^{k,t}_{x_i}   & \text{if $\alpha^{k,t}_{c,x_i} \leq \epsilon_c$ and $\beta^{k,t}_{c,x_i} > \epsilon_c$} \\
\hat \rho_{x_i}^{k,t}   & \text{if $\alpha^{k,t}_{c,x_i} \leq \epsilon_c$ and $\beta^{k,t}_{c,x_i} \leq \epsilon_c$}
\end{array} \right.   \nonumber \\
 \rho_{x_i}^{k,t+1}  &=&  \max \{\min\{\hat \rho_{x_i}^{k,t+1}, U_{x_i}\}, L_{x_i}\}, \ \ \ x \in \{v,\theta\},
\end{eqnarray}
where $U_{x_i}$ and $L_{x_i}$ are the upper and lower bounds of the penalty parameter. We can similarly calculate $\rho_{p_e}$ and $\rho_{q_e}$ for every link $e$ with $|\mathcal M_e| > 1$.
\section{Numerical Results} \label{sec:num}
An open-source software package DiCARP \cite{Karimi_DiCARP} is created for the algorithms in this paper. The algorithms are implemented in Python using the optimization modeling package Pyomo \cite{hart2011pyomo,bynum2021pyomo}, with the backend solver Ipopt \cite{ipopt}. Ipopt is an open source software package for large-scale nonlinear optimization. Our DiCARP package contains a method that applies the greedy radial partitioning to a given graph. The problem instances we solve are from the MATPOWER library \cite{zimmerman2010matpower}. A function is given to transform a MATPOWER case into Pyomo's data file format. We performed computational experiments on a 1.7 GHz 12th Gen Intel Core i7 personal computer with 32GB of memory. In this section, we refer to the solution returned by Ipopt as $P_{IPM}$ and the solution returned by the DiCA as $P_{DiCA}$. To compare these two quantities, similar to \cite{mhanna2018component}, we define an optimality $GAP$ as:
\begin{eqnarray} \label{eqn:gap}
GAP = \left|{P_{IPM} - P_{DiCA}} \right|/{P_{IPM}}
\end{eqnarray}
For updating the penalty parameters in Algorithm \ref{alg:DiCA}, we use the Spectral Parameter Selection described in Section \ref{sec:adap}. Since the approach is adaptive, the sensitivity of the algorithm to the penalty parameters' initial values are lower compared to non-adaptive DO. For numerical results, we set $\rho^{k,0}_{v_i}=\rho^{k,0}_{\theta_i}=10^4$ for all the nodes with $|\mathcal M_i| >1$, and $\rho^{k,0}_{p_e}=\rho^{k,0}_{q_e}=10^3$ for all the edges with $|\mathcal M_e| > 1$. 

The most comprehensive numerical results for component-based DO for power networks are \cite{mhanna2018component} for non-adaptive case and \cite{mhanna2018adaptive} for adaptive DO. Table \ref{tab:res1} shows the results of solving problems from the MATPOWER library that were also reported in \cite{mhanna2018component}. As can be seen, the group-based DiCA returns more accurate solutions in much smaller number of iterations, without the requirement of tuning the parameters.

\begin{table}[!t]
\caption{Performance of the DiCA algorithm on MATPOWER instances \label{tab:res1}}
\centering
\begin{tabular}{|c||c|c|c|c||c|c|}
\hline
Case & iter  & $P_{DiCA}$ & $P_{IPM}$ & GAP  & iter \cite{mhanna2018component}  & GAP \cite{mhanna2018component} \\
\hline \hline
5 & 248 &  17551.89 & 17551.89 & 4.51 e-9    & 3911 & 4.14e-5 \\
\hline
6ww & 64 &  3143.97 & 3143.97 & 2.12e-8  & 918  & 4.94e-5 \\
\hline
9  &  44 &  5296.68 & 5296.68 & 1.13e-08 & 630 &  7.25e-5 \\
\hline
14 & 72 &  8081.52 & 8081.52 & 3.53e-08 & 857 &  7.88e-5  \\
\hline
24 & 115 &   63352.20 &  63352.20 & 2.38e-08  & 924 & 2.05e-06  \\
\hline
30 & 532 &  576.89 & 576.89 & 7.74e-07 & 2763 & 1.34e-4  \\  \hline
39 & 342 &  41864.14 & 41864.18 & 1.28 e-08  & 7468  & 9.82e-07\\
\hline 
57 & 232 &   41737.79 & 41737.78 & 2.39e-07 & 1305 & 4.24e-05 \\ \hline
118 & 215  &  129660.81 & 129660.69 & 9.25e-07 & 1168  & 3.66e-6  \\ \hline
300 & 684 &   719725.54 & 719725.09 & 6.25e-07 & 11755 & 8.72e-7 \\ \hline
\end{tabular}
\end{table}

 Table \ref{tab:res2} reports the results of using DiCA for larger instances that are solved using the adaptive component-based DO in \cite{mhanna2018adaptive}. The DiCA algorithm performs better compared to the adaptive component-based algorithm. However, the performance gap gets smaller as the radial partitioning returns more number of partitions.  Our numerical results align with the limited literature on the topic of partitioning, indicating that reducing the number of sub-problems can significantly decrease the number of iterations. We did not report the running times in our tables, since IPM is much faster using a single machine for experiments, and the running times are not informative. Using distributed optimization (DO) is justified when data privacy necessitates multiple agents or the problem size is so large that a generic solver like Ipopt cannot handle it. These are the topics of future research for very large modern power networks. 
 
\begin{table}[!t]
\caption{Performance of the DiCA algorithm on MATPOWER instances also solved in \cite{mhanna2018adaptive} using adaptive component-based DO. \label{tab:res2}}
\centering
\begin{tabular}{|c||c|c||c|c|}
\hline
Case & iter  &  GAP  & iter \cite{mhanna2018adaptive}  & GAP \cite{mhanna2018adaptive} \\
\hline 
1354 \_PEGASE & 753 &  6.75 e-7    &  1110 & 6.43e-6 \\
\hline
2383\_wp & 1740 &   7.81e-7  & 4070  & 1.95e-5 \\ \hline
2736\_sp & 1212 & 5.42e-7  & 2154  & 1.33e-6 \\ \hline
2746\_wp & 986 &   3.21e-6  & 1872  & 5.09e-5 \\ 
\hline
\end{tabular}
\end{table}

\begin{figure}[htp]
    \centering
    \begin{tabular}{cc}
        \begin{subfigure}[b]{0.45\textwidth}
            \centering
            \includegraphics[width=\textwidth]{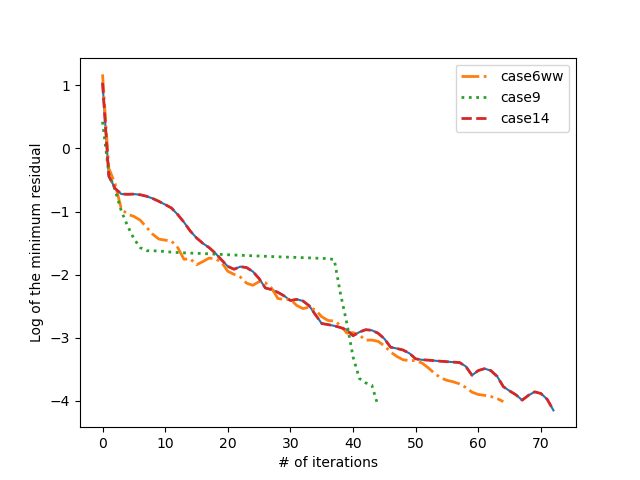}
            \caption{}
            \label{fig:fig1}
        \end{subfigure} &
        \begin{subfigure}[b]{0.45\textwidth}
            \centering
            \includegraphics[width=\textwidth]{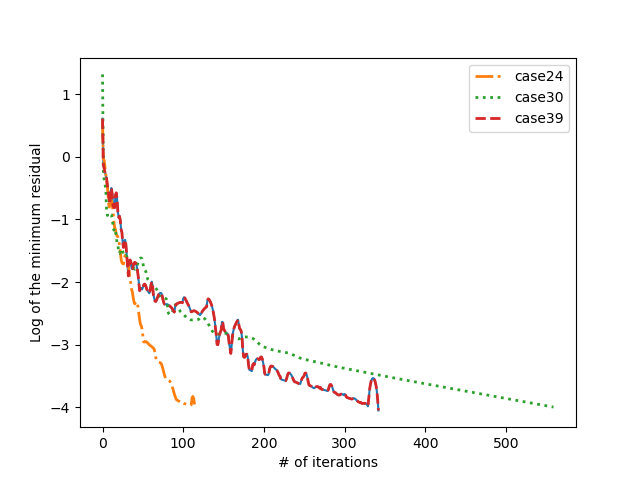}
            \caption{}
            \label{fig:fig2}
        \end{subfigure} \\
        \begin{subfigure}[b]{0.45\textwidth}
            \centering
            \includegraphics[width=\textwidth]{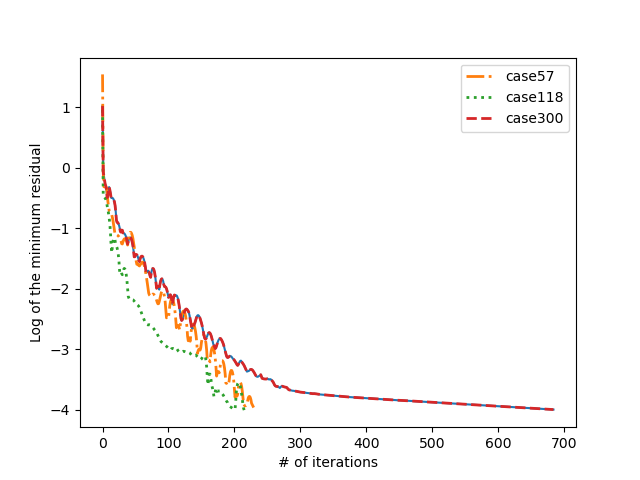}
            \caption{}
            \label{fig:fig3}
        \end{subfigure} &
        \begin{minipage}[b]{0.45\textwidth}
            \centering
    
            \caption{The minimum residual of the regions versus the number of iterations for different problems from the MATPOWER library (a) case6ww, case9 and case14, (b) case24, case30, and case39, and (c) case57, case118, and case300.}
            \label{fig:res_rho}
            \vspace{10pt}
        \end{minipage}
    \end{tabular}
\end{figure}

Figure \ref{fig:res_rho} shows the plots of the minimum residual of the regions versus the number of iterations. As explained in \cite{mhanna2018component}, the oscillations in the progress of the residual is one of the reasons of slow convergence in DO. As can be seen, for almost all examples, we do not see oscillations in the plots performing our DiCA algorithm.   
%
%

One of the benefits of using radial partitioning is simplifying the structure of the sub-problems. One way to measure this is by counting the number of iterations for solving the problem using a generic interior-point solver. This works since in practice, the number of iterations of interior-point methods is independent of the size of the problem, and relies more on the structure \cite{nocedal1999numerical, boyd2004convex}.  In other words, problems that are closer to be \emph{ill-conditioned} take more iterations to solve by interior-point methods.  Table \ref{tab:res3} shows that number of iterations for solving some problem instances using Ipopt, and the average number of iterations for the sub-problems, which are significantly  smaller.

\begin{table}[!t]
\caption{The number of iterations of Ipopt  \label{tab:res3}}
\centering
\begin{tabular}{|c||c|c|c|c|c|c|}
\hline
Case &  9 & 14& 30 & 57 & 118 & 300 \\
\hline \hline
Main problem & 30 &  48 & 69 & 53  & 88 & 138 \\ \hline
Average of sub-problems &  25.5 & 25.3 & 23.8  & 25.6 & 30.3 & 33.4  \\
\hline

\end{tabular}
\end{table}
\section{Conclusion} \label{sec:con}
In this work, we proposed an approach in partitioning to reduce the number of sub-problems and enforce desirable structures to them. The numerical results show that our DiCA algorithm applied on top of our proposed radial partitioning and with spectral parameter selection returns more accurate solutions in fewer iterations than the adaptive component-based algorithms. The spectral parameter selection alleviate the need complicated parameter tuning. 

Many open questions arise in the context of radial and, in general, intelligent partitioning. We observed that, on average, the sub-problems take less than half of the number of iterations of the main problem using a generic interior-point solver. This fact is evidence that radial partitioning simplifies the structure of the sub-problems. The next step is to prove the existence of exact convex relaxations for the sub-problems, using the current results on convex relaxations for radial networks, such as \cite{gan2014exact,wei2017optimal,farivar2013branch}. Another future work is designing radial partitionings other than our proposed greedy algorithm that consider the power network structure more efficiently. Another issue shown by our numerical results is the number of partitions created by the radial partitioning is increasing by the size of the problem, which brings back the slow convergence for very large problems. An interesting open question is finding an intelligent partitioning that returns very small number of partitions even for large problem instances. For parameter selection, an interesting next step is studying different adaptive or ML-based parameter selection methods. 
 
 \section{Acknowledgements}
 
 Research of the author was supported by	 the US National
Science Foundation (NSF) under Grant No. CMMI-2347120.


  \bibliographystyle{siam} 
  \bibliography{references}

\end{document}